\documentclass{amsart}
\usepackage[dvips]{graphicx}
\usepackage{amscd}
\usepackage{amsmath}
\usepackage{amsxtra}
\usepackage{amsfonts}
\usepackage{amssymb}

\newtheorem{theorem}{Theorem}[section]
\newtheorem{corollary}[theorem]{Corollary}
\newtheorem{lemma}[theorem]{Lemma}

\theoremstyle{definition}
\newtheorem{definition}[theorem]{Definition}

\newtheorem{example}[theorem]{Example}
\theoremstyle{remark}

\renewcommand{\theclaim}{\textup{\theclaim}}

\newtheorem*{acknowledgements}{Acknowledgements}

\renewcommand{\theenumi}{\roman{enumi}}

\numberwithin{equation}{section}

\def\hooklongrightarrow{\lhook\joinrel\longrightarrow} 

\def\openone
{\mathchoice
{\hbox{\upshape \small1\kern-3.3pt\normalsize1}}
{\hbox{\upshape \small1\kern-3.3pt\normalsize1}}
{\hbox{\upshape \tiny1\kern-2.3pt\SMALL1}}
{\hbox{\upshape \Tiny1\kern-2pt\tiny1}}}
\makeatletter
\newbox\ipbox
\newcommand{\ip}[2]{\left\langle #1\mathrel{\mathchoice
{\setbox\ipbox=\hbox{$\displaystyle \left\langle\mathstrut #1#2\right\rangle$}
\vrule height\ht\ipbox width0.25pt depth\dp\ipbox}
{\setbox\ipbox=\hbox{$\textstyle \left\langle\mathstrut #1#2\right\rangle$}
\vrule height\ht\ipbox width0.25pt depth\dp\ipbox}
{\setbox\ipbox=\hbox{$\scriptstyle \left\langle\mathstrut #1#2\right\rangle$}
\vrule height\ht\ipbox width0.25pt depth\dp\ipbox}
{\setbox\ipbox=\hbox{$\scriptscriptstyle \left\langle\mathstrut #1#2\right\rangle$}
\vrule height\ht\ipbox width0.25pt depth\dp\ipbox}
} #2\right\rangle}
\newcommand{\diracb}[1]{\left\langle #1\mathrel{\mathchoice
{\setbox\ipbox=\hbox{$\displaystyle \left\langle\mathstrut #1\right.$}
\vrule height\ht\ipbox width0.25pt depth\dp\ipbox}
{\setbox\ipbox=\hbox{$\textstyle \left\langle\mathstrut #1\right.$}
\vrule height\ht\ipbox width0.25pt depth\dp\ipbox}
{\setbox\ipbox=\hbox{$\scriptstyle \left\langle\mathstrut #1\right.$}
\vrule height\ht\ipbox width0.25pt depth\dp\ipbox}
{\setbox\ipbox=\hbox{$\scriptscriptstyle \left\langle\mathstrut #1\right.$}
\vrule height\ht\ipbox width0.25pt depth\dp\ipbox}
}\right. }
\newcommand{\dirack}[1]{\left. \mathrel{\mathchoice
{\setbox\ipbox=\hbox{$\displaystyle \left.\mathstrut #1\right\rangle$}
\vrule height\ht\ipbox width0.25pt depth\dp\ipbox}
{\setbox\ipbox=\hbox{$\textstyle \left.\mathstrut #1\right\rangle$}
\vrule height\ht\ipbox width0.25pt depth\dp\ipbox}
{\setbox\ipbox=\hbox{$\scriptstyle \left.\mathstrut #1\right\rangle$}
\vrule height\ht\ipbox width0.25pt depth\dp\ipbox}
{\setbox\ipbox=\hbox{$\scriptscriptstyle \left.\mathstrut #1\right\rangle$}
\vrule height\ht\ipbox width0.25pt depth\dp\ipbox}
} #1\right\rangle}
\newcommand{\lonet}{L^{1}\left(  \mathbb{T}\right)}

\newcommand{\linft}{L^{\infty}\left(  \mathbb{T}\right)}

\newcommand{\ltwor}{L^{2}\left(\mathbb{R}\right)}
\newcommand{\linfr}{L^{\infty}\left(\mathbb{R}\right)}
\input cyracc.def

\begin{document}
\title[Signed Ruelle operators]{Harmonic analysis of signed Ruelle transfer operators}
\author{Dorin Ervin Dutkay}
\address{Department of Mathematics\\
The University of Iowa\\
14 MacLean Hall\\
Iowa City, IA 52242-1419\\
U.S.A.}
\email{ddutkay@math.uiowa.edu}
\thanks{}
\subjclass{}
\keywords{}

\begin{abstract}
Motivated by wavelet analysis, we prove that there is a one-to-one correspondence between the following data:
\begin{enumerate}
\item 
Solutions to $R(h)=h$ where $R$ is a certain non-positive Ruelle transfer operator;
\item
Operators that intertwine a certain class of representations of the $C^*$-algebra $\mathfrak{A}_N$ 
on two unitary generators $U$, $V$ subject to the relation 
$$UVU^{-1}=V^N$$ 
\end{enumerate}
This correspondence enables us to give a criterion for the biorthogonality of 
a pair of scaling functions and calculate all solutions of the equation $R(h)=h$ in some concrete cases. 
\end{abstract}\maketitle
\tableofcontents
\section{\label{Intr}Introduction}
The multiresolution wavelet theory establishes a close interconnection 
between two operators: $M$ - the cascade refinement operator and $R$ - the transfer operator,
also called Ruelle operator ( see \cite{Dau92}, \cite{Jor} ). Our present approach stresses 
representation theory and intertwining operators. 
\par
In this paper we show how to get wavelets from representations and we compare representations which yield 
different wavelets. Examples are given in section 4. 
\par
We recall that $M$ operates on $L^2\left(\mathbb{R}\right)$ by
$$M\psi(x)=\sqrt{N}\sum_{k\in\mathbb{Z}}a_k\psi(Nx-k),\quad(x\in\mathbb{R})$$
or, equivalently, in Fourier space 
$$\widehat{M\psi}(x)=\frac{m_0\left(\frac{x}N\right)}{\sqrt{N}}\widehat\psi\left(\frac{x}N\right),\quad (x\in\mathbb{R})$$
where $N\geq 2$ is an integer - the scale, $m_0(z)=\sum_{k\in\mathbb{Z}}a_kz^k$ for $z\in\mathbb{T}$, $\mathbb{T}$ being the unit circle, and $\widehat{\psi}$ denotes 
the Fourier transform 
$$\widehat{\psi}(x)=\int_{\mathbb{R}}\psi(t)e^{-itx}\,dt.$$

\par
The Ruelle transfer operator is defined on $\lonet$ by
$$Rf(z)=\frac{1}{N}\sum_{w^N=z}\left|m_0(w)\right|^2f(w),\quad(z\in\mathbb{T}).$$
On $\mathbb{T}$, we consider $\mu$, the normalized Haar measure.  
\par
It is the equation
\begin{equation}\label{eqscale}
M\varphi=\varphi,
\end{equation}
or, equivalently,
$$\varphi(x)=\sqrt{N}\sum_{k\in\mathbb{Z}}a_k\varphi(Nx-k),\quad(x\in\mathbb{R})$$
which generates the wavelets. It is called the refinement (or scaling) equation. 
\par
 The orthogonality properties of the integer translates of the scaling function $\varphi\in L^2\left(\mathbb{R}\right)$, $M\varphi=\varphi$ 
are closely connected to the problem of finding a positive eigenvector for $R$ 
\begin{equation}\label{eigen}
h\in\lonet\, ,\, h\geq0\, ,\, Rh=h
\end{equation} 
 (see \cite{CoDa92}, \cite{BrJo99}, \cite{CoRa90} where a correspondence is established between the non-zero $\ltwor$-solutions $\varphi$ 
to (\ref{eqscale}) and the non-zero solutions $h$ to (\ref{eigen}). In general, solutions need not exist.) 
A necessary condition for the orthogonality of the translates of the scaling function is the quadrature mirror 
filter restriction: 
$$\frac{1}{N}\sum_{w^N=z}\left|m_0(w)\right|^2=1\quad (z\in\mathbb{T})$$
which , in terms of the Ruelle operator can be rewritten as:
$$R\openone=\openone.$$
\par
Lawton( \cite{Law91a}) gave a neccesary and sufficient condition formulated also in terms of the Ruelle operator: 
the translates of the scaling function are orthogonal if and only if the constant function $\openone$ is the only 
continuous solution of (\ref{eigen}) (up to a multiplicative costant). 
\par
The scaling equation (\ref{eqscale}) can be reformulated in a $C^*$-algebra setting.
\par
Consider $\mathfrak{A}_N$, the $C^*$-algebra generated by two unitary operators $U$ and $V$, 
satisfying the relation $UVU^{-1}=V^N$. It has a representation on 
$L^2\left(\mathbb{R}\right)$ given by 
$$U :\psi\mapsto\frac{1}{\sqrt{N}}\psi\left(\frac{x}{N}\right)\quad ,\quad 
V:\psi\mapsto\psi(x-1)\quad(x\in\mathbb{R})$$
$V=\pi(z)$ where $\pi$ is the representation of $\linft$ given by
$$\left(\pi(f)\psi\right)^{\widehat{ }}=f\widehat\psi\, , \quad (f\in\linft).$$
\par
The scaling equation (\ref{eqscale}) becomes 
$$U\varphi=\pi\left(m_0\right)\varphi.$$
\par
The system $(U,\pi,L^2\left(\mathbb{R}\right),\varphi,m_0)$ is called 
the wavelet representation with scaling function $\varphi$ (see \cite{Jor}).
\par
If a wavelet representation is given with scaling function $\varphi$ then it produces a solution for (\ref{eigen}):
$$h_{\varphi}(z)=\sum_{n\in\mathbb{Z}}z^n\ip{\pi\left(z^n\right)\varphi}{\varphi}=
\frac{1}{2\pi}\sum_{k\in\mathbb{Z}}\left|\widehat{\varphi}(\omega+2k\pi)\right|^2,\quad(z=e^{-i\omega}).$$

In \cite{Jor} it is proved that a converse also holds. Any solution $h\geq0$ to $Rh=h$ arises in this way, as 
$h=h_{\varphi}$ for some representation $\pi$ of $\mathfrak{A}_N$. 
\par
Thus, the analysis of orthogonal wavelets is closely related to the study of the positive Ruelle operator $R$ and this operator 
is linked to the representations of the algebra $\mathfrak{A}_N$. 
\par
For an analysis of biorthogonal wavelets, it turns out that we have to consider non-positive Ruelle operators. They 
correspond to a pair of filters $m_0,m_0'\in\linft$ and are defined by:
$$R_{m_0,m_0'}f(z)=\frac{1}{N}\sum_{w^N=z}\overline{m_0(w)}m_0'(w)f(w),\quad (f\in\lonet\, ,\, z\in\mathbb{T} ).$$ 
\par
The condition corresponding to the quadrature mirror filter condition, and necessary for the biorthogonality of wavelets, is
$$\frac{1}{N}\sum_{w^N=z}\overline{m_0(w)}m_0'(w)=1,\quad (z\in\mathbb{T})$$
which rewrites as
$$R_{m_0,m_0'}\openone=\openone.$$
\par
If two scaling functions $\varphi,\varphi'$ are given, with $U\varphi=\pi\left(m_0\right)\varphi$, 
$U\varphi'=\pi\left(m_0'\right)\varphi'$, then 
$$h_{\varphi,\varphi'}(z)=\sum_{n\in\mathbb{Z}}z^n\ip{\pi\left(z^n\right)\varphi}{\varphi'}=
\frac{1}{2\pi}\sum_{k\in\mathbb{Z}}\overline{\widehat{\varphi}}\widehat{\varphi'}(\omega+2k\pi),\quad(z=e^{-i\omega})$$
satisfies 
$$R_{m_0,m_0'}h_{\varphi,\varphi'}=h_{\varphi,\varphi'}.$$
\par For more background on wavelets we refer the reader to \cite{Dau92}. 
\par
We will see in this paper that solutions to $R_{m_0,m_0'}h=h$ correspond to operators 
that intertwine the representations of $\mathfrak{A}_N$ introduced in \cite{Jor} arrising from $m_0$ and $m_0'$, respectively.
 In chapter \ref{Main} we establish this correspondence and in chapter \ref{Applications} 
we give a criterion for the biorthogonality of two given scaling functions in terms of 
the eigenspace of the non-positive Ruelle transfer operator $R_{m_0,m_0'}$ associated to the eigenvalue $1$. In chapter 
\ref{Examples} we consider some concrete examples of filters and give complete solutions for 
the equation $Rh=h$.

\section{\label{Main} Main results}
In this section we prove our main theorems on wavelets and representations: theorem \ref{th1} and theorem \ref{th3}.
These results prove the bijective correspondence between two sets: operators that intertwine the cyclic representations presented in 
\cite{Jor} and solutions to $R_{m_0,m_0'}h=h$ .
\par 
We begin with some properties of the Ruelle operator. We will denote by $R=R_{m_0,m_0'}$, 
$m_0,m_0'\in\linft$. 
\begin{lemma} 
For $f\in\lonet$ 
\label{Elemprop}
\begin{enumerate}
\item \label{Elemprop1}
$$\int_{\mathbb{T}}Rf(z)\,d\mu=\int_{\mathbb{T}}\overline{m_0(z)}m_0'(z)f(z)\,d\mu.$$
\item \label{Elemprop2}
$$\int_{\mathbb{T}}g(z)Rf(z)\,d\mu=\int_{\mathbb{T}}g(z^N)\overline{m_0(z)}m_0'(z)f(z)\,d\mu.$$
\item \label{Elemprop3}
$$R(g(z^N)f(z))=g(z)Rf(z),$$
$$R^n(g(z^{N^n})f(z))=g(z)R^nf(z).$$
\item \label{Elemprop4}
$$\int_{\mathbb{T}}R^nf(z)\,d\mu=\int_{\mathbb{T}}\overline{m_0^{(n)}(z)}{m_0'}^{(n)}f(z)\,d\mu.$$
where $m_0^{(n)}(z)=m_0(z)m_0\left(z^N\right)\ldots m_0\left(z^{N^{n-1}}\right)$ .
\end{enumerate}
\end{lemma}
\begin{proof}

(i)
\begin{align*}
\int_{\mathbb{T}}Rf(z)\,d\mu & =\frac 1N\sum_{k=0}^{N-1}\frac{1}{2\pi}\int_0^{2\pi}\overline{m_0\left(\frac{\theta+2k\pi}{N}\right)}m_0'\left(\frac{\theta+2k\pi}{N}\right)f\left(\frac{\theta+2k\pi}{N}\right)\,d\theta\\
 & =\sum_{k=0}^{N-1}\frac{1}{2\pi}\int_{\frac{2k\pi}{N}}^{\frac{2(k+1)\pi}N}\overline{m_0(\theta)}m_0'(\theta)f(\theta)\,d\theta\\
& =\int_{\mathbb{T}}\overline{m_0(z)}m_0'(z)f(z)\,d\mu.
\end{align*}
(iii) Clear.\\
(ii) Follows from (i) and (iii).\\
(iv) Proof by induction. For $n=1$ it is (i). 
\begin{align*}
\int_{\mathbb{T}}R^{n+1}f\,d\mu & =\int_{\mathbb{T}}R\left(R^nf\right)\,d\mu=\int_{\mathbb{T}}\overline{m_0(z)}m_0'(z)R^nf(z)\,d\mu\\
& =\int_{\mathbb{T}}R^n\left(\overline{m_0\left(z^{N^n}\right)}m_0'\left(z^{N^n}\right)f(z)\right)\,d\mu\\
& = \int_{\mathbb{T}}\overline{m_0^{(n)}(z)}{m_0'}^{(n)}(z)\overline{m_0\left(z^{N^n}\right)}m_0'\left(z^{N^n}\right)f(z)\,d\mu\\
& = \int_{\mathbb{T}}\overline{m_0^{(n+1)}(z)}{m_0'}^{(n+1)}(z)f(z)\,d\mu.
\end{align*}
\end{proof}
From \cite{Jor} theorem 2.4 we know that, given $m_0\in\linft$ which is non-singular (i.e. doesn't vanish on a 
subset of positive measure), there is a 1-1 correspondence between 
$$(a)\,\,\,\, h\in\lonet\, , \, h\geq 0\, , \, R(h)=h\,\,\,\, (\mbox{here }R=R_{m_0,m_0})$$
and
$$(b) \tilde{\pi}\in\operatorname*{Rep}\left(\mathfrak{A}_N,\mathcal{H}\right), \, \varphi\in\mathcal{H}$$
with the unitary $U$ from $\tilde{\pi}$ satisfying 
$$U\varphi=\pi\left(m_0\right)\varphi$$ 
$\operatorname*{Rep}\left(\mathfrak{A}_N,\mathcal{H}\right)$ 
is the set of normal representations of the algebra $\mathfrak{A}_N$. These representations are in fact generated by 
a unitary $U$ on $\mathcal{H}$ and a representation $\pi$ of $\linft$ on $\mathcal{H}$, with the property that 
$$U\pi(f(z))U^{-1}=\pi\left(f\left(z^N\right)\right),\quad(f\in\linft).$$

Here is again theorem 2.4 from \cite{Jor}:
\begin{theorem}
\label{Thmax+b.3}\ 

\begin{enumerate}
\item \label{Thmax+b.3(1)}Let $m_{0}\in L^{\infty}\left(  \mathbb{T}\right)
$, and suppose $m_{0}$ does not vanish on a subset of $\mathbb{T}$ of positive
measure. Let%
\begin{equation}
\left(  Rf\right)  \left(  z\right)  =\frac{1}{N}\sum_{w^{N}=z}\left|
m_{0}\left(  w\right)  \right|  ^{2}f\left(  w\right)  ,\qquad f\in
L^{1}\left(  \mathbb{T}\right)  . \label{eqax+b.14}%
\end{equation}
Then there is a one-to-one correspondence between the data
\textup{(\ref{Thmax+b.3(1)(1)})} and \textup{(\ref{Thmax+b.3(1)(2)})} below,
where \textup{(\ref{Thmax+b.3(1)(2)})} is understood as equivalence classes
under unitary equivalence:

\begin{enumerate}
\renewcommand{\theenumi}{\relax}

\item \label{Thmax+b.3(1)(1)}$h\in L^{1}\left(  \mathbb{T}\right)  $, $h\geq
0$, and%
\begin{equation}
R\left(  h\right)  =h. \label{eqax+b.15}%
\end{equation}

\item \label{Thmax+b.3(1)(2)}$\tilde{\pi}\in\operatorname*{Rep}\left(
\mathfrak{A}_{N},\mathcal{H}\right)  $, $\varphi\in\mathcal{H}$, and the unitary
$U$ from $\tilde{\pi}$ satisfying%
\begin{equation}
U\varphi=\pi\left(  m_{0}\right)  \varphi. \label{eqax+b.16}%
\end{equation}
\end{enumerate}

\item \label{Thmax+b.3(2)}From \textup{(\ref{Thmax+b.3(1)(1)})}$\rightarrow
$\textup{(\ref{Thmax+b.3(1)(2)}),} the correspondence is given by%
\begin{equation}
\ip{\varphi}{\pi\left( f\right) \varphi}_{\mathcal{H}}=\int_{\mathbb{T}%
}fh\,d\mu, \label{eqax+b.17}%
\end{equation}
where $\mu$ denotes the normalized Haar measure on $\mathbb{T}$.

 From \textup{(\ref{Thmax+b.3(1)(2)})}$\rightarrow$%
\textup{(\ref{Thmax+b.3(1)(1)}),} the correspondence is given by%
\begin{equation}
h\left(  z\right)  =h_{\varphi}\left(  z\right)  =\sum_{n\in\mathbb{Z}}%
z^{n}\ip{\pi\left( e_{n}\right) \varphi}{\varphi}_{\mathcal{H}}.
\label{eqax+b.18}%
\end{equation}

\item \label{Thmax+b.3(3)}When \textup{(\ref{Thmax+b.3(1)(1)})} is given to
hold for some $h$, and $\tilde{\pi}\in\operatorname*{Rep}\left(  \mathfrak{A}%
_{N},\mathcal{H}\right)  $ is the corresponding cyclic representation with
$U\varphi=\pi\left(  m_{0}\right)  \varphi$, then the representation is unique
from $h$ and \textup{(\ref{eqax+b.17})} up to unitary equivalence: that is, if
$\tilde{\pi}^{\prime}\in\operatorname*{Rep}\left(  \mathfrak{A}_{N},\mathcal{H}^{\prime
}\right)  $, $\varphi^{\prime}\in\mathcal{H}^{\prime}$ also cyclic and
satisfying%
\begin{align*}
\ip{\varphi^{\prime}}{\pi^{\prime}\left( f\right) \varphi^{\prime}} &
=\int_{\mathbb{T}}fh\,d\mu\\%
\intertext{and}%
U^{\prime}\varphi^{\prime} &  =\pi^{\prime}\left(  m_{0}\right)
\varphi^{\prime},
\end{align*}
then there is a unitary isomorphism $W$ of $\mathcal{H}$ onto $\mathcal{H}%
^{\prime}$ such that $W\pi\left(  f\right)  =\pi^{\prime}\left(  f\right)  W$, for
$f\in\linft$, $WU=U^{\prime}W$ and $W\varphi=\varphi^{\prime}$.
\end{enumerate}
\end{theorem}
\begin{definition}
\label{cyclicreph}
Given $h$ as in theorem \ref{Thmax+b.3} call $\left(\pi, U, \mathcal{H}, \varphi\right)$ the cyclic representation of $\mathfrak{A}_N$ 
associated to $h$.
\end{definition}
The next theorem shows how solutions of $R_{m_0,m_0'}h_0=h_0$ induce operators that intertwine these cyclic representations.
\begin{theorem}
\label{th1}
Let $m_0,m_0'\in\linft$ be non-singular and $h,h'\in\lonet$, $h,h'\geq0$, $R_{m_0,m_0}(h)=h$,
$R_{m_0',m_0'}(h')=h'$. Let $\left(\pi,U,\mathcal{H},\varphi\right)$, $\left(\pi ',U',\mathcal{H}',\varphi '\right)$
be the cyclic representations corresponding to $h$ and $h'$ respectively. \\
If $h_0\in\lonet$, $R_{m_0,m_0'}\left(h_0\right)=h_0$ and $\left| h_0\right| ^2\leq chh'$ for some $c>0$ then 
there exists a unique operator $S:\mathcal{H}'\rightarrow\mathcal{H}$ such that 
$$SU'=US, \quad S\pi '(f)=\pi (f)S, \quad (f\in\linft)$$
$$\ip{\varphi}{\pi (f)S\varphi '}=\int_{\mathbb{T}}fh_0\,d\mu,\quad(f\in\linft)$$
Moreover $\left\| S\right\|\leq\sqrt{c}$.
\end{theorem}
\begin{proof}
To symplify the notation let $R_0:=R_{m_0,m_0'}$. Look at the construction of $\tilde{\pi}$ and $\mathcal{H}$ 
in the proof of theorem 2.4. in \cite{Jor}. We reproduce here the main steps of this construction. 
First, one considers 
\[
\mathcal{V}_{n}:=\left\{  \left(  \xi,n\right)  \mid\xi\in L^{\infty}\left(
\mathbb{T}\right)  \right\}
\]
and
$$\ip{\left( \xi,n\right) }{\left( \eta,n\right) }_{\mathcal{H}}  
=\int_{\mathbb{T}}R^{n}\left(  \bar{\xi}\eta h\right)  \,d\mu\text{\qquad for
}n=1,2,\dots $$
Let $\mathcal{H}_n$ be the completion of $\mathcal{V}_n$ in this scalar product. 
\par
When $n$, $k$ are given, $n\geq0$, $k\geq1$, one constructs the isometry
$\mathcal{V}_{n}\hookrightarrow\mathcal{V}_{n+k}$ by iteration of the one from
$\mathcal{V}_{n}$ to $\mathcal{V}_{n+1}$, i.e.,%
\[
\mathcal{V}_{n}\hooklongrightarrow\mathcal{V}_{n+1}\hooklongrightarrow
\mathcal{V}_{n+2}\hooklongrightarrow\cdots\hooklongrightarrow\mathcal{V}%
_{n+k},
\]
where $J\colon\mathcal{V}_{n}\rightarrow\mathcal{V}_{n+1}$ is defined by
$$
J\left(  \left(  \xi,n\right)  \right)  :=\left(  \xi\left(  z^{N}\right)
,n+1\right)  . \label{eqPoof.9}%
$$
Then $\mathcal{H}$ is defined as the inductive limit of the Hilbert spaces $\mathcal{H}_n$. The set 
$\cup_{n\geq0}\mathcal{V}_n$ is dense in $\mathcal{H}$. 
\par
The representation is defined as follows 

$$U\left(  \xi,0\right)     :=\left(  S_{0}\xi,0\right)  =\left(  m_{0}\left(
z\right)  \xi\left(  z^{N}\right)  ,0\right)  ,$$
$$U\left(  \xi,n+1\right)     :=\left(  m_{0}\left(  z^{N^{n}}\right)
\xi\left(  z\right)  ,n\right)  , $$
and
$$\pi\left(  f\right)  \left(  \xi,n\right)     :=\left(  f\left(  z^{N^{n}}\right)  \xi\left(  z\right)  ,n\right)  
$$
\par
The scaling function $\varphi$ is $(\openone,0)$. Recall also the main property of this representation 
$$U\pi(f)=\pi\left(f\left(z^N\right)\right)U,\quad (f\in\linft).$$
\par
Having these, we return to our proof.
Define
$$B\left[(\xi , n)|(\xi ',n)\right]=\int_{\mathbb{T}}R_0^n\left(\overline{\xi}\xi 'h_0\right)\,d\mu\, ,\quad\mbox{for}\, 
(\xi,n)\in\mathcal{V}_n \, ,\, (\xi ',n)\in\mathcal{V}_n'.$$
Then
\begin{align*}
\left| B\left[(\xi , n)|(\xi ' ,n)\right]\right| &= \left|\int_{\mathbb{T}}R_0^n\left(\overline{\xi}\xi 'h_0\right)\,d\mu\right| \\
&\leq\int_{\mathbb{T}}\left|\overline{m_0^{(n)}(z)}{m_0'}^{(n)}(z)\overline{\xi}\xi 'h_0\right|\,d\mu\quad\mbox{(by (\ref{Elemprop4}) of lemma \ref{Elemprop})}\\
&\leq\sqrt{c}\int_{\mathbb{T}}\left|m_0^{(n)}(z)\right|\left|{m_0'}^{(n)}(z)\right|\left|\xi '\right|\left|\xi\right|\left| h\right|^{\frac 12}\left|h'\right|^{\frac 12}\,d\mu\\
&\leq\sqrt{c}\left(\int_{\mathbb{T}}\left|m_0^{(n)}(z)\right|^2\left|\xi\right|^2h\,d\mu\right)^{\frac 12}
\left(\int_{\mathbb{T}}\left|{m_0'}^{(n)}\right|^2\left|\xi '\right|^2h'\,d\mu\right)^{\frac 12}\\
&=\sqrt{c}\left\|(\xi , n)\right\|_{\mathcal{H}}\left\|(\xi ',n)\right\|_{\mathcal{H}'}.
\end{align*}
Therefore
\begin{equation}
\left| B\left[ (\xi , n)|(\xi ',n)\right]\right|\leq\sqrt{c}\left\| (\xi , n)\right\|_{\mathcal{H}}\left\| (\xi ',n)\right\|_{\mathcal{H}'}.
\label{eq1}
\end{equation}
Equation (\ref{eq1}) implies also that $B$ can be extended from $\mathcal{V}_n\times\mathcal{V}_n'$ to 
$\mathcal{H}_n\times\mathcal{H}_n'$ such that (\ref{eq1}) remains valid. 
\par 
Next we prove that $B$ is compatible with the inductive limit structure that define the Hilbert spaces $\mathcal{H}$
 and $\mathcal{H}'$.
\begin{align*}
B\left[J(\xi , n)|J'(\xi ',n)\right]&=B\left[\left(\xi (z^N),n+1\right)|\left(\xi '(z^N),n+1\right)\right]\\
&=\int_{\mathbb{T}}R_0^{n+1}\left(\overline{\xi (z^N)}\xi '(z^N)h_0(z)\right)\,d\mu\\
&=\int_{\mathbb{T}}R_0^n\left(R_0\left(\overline{\xi (z^N)}\xi '(z^N)h_0(z)\right)\right)\,d\mu\\
&=\int_{\mathbb{T}}R_0^n\left(\overline{\xi (z)}\xi '(z)R_0\left( h_0\right)\right)\,d\mu\\
&=\int_{\mathbb{T}}R_0^n\left(\overline{\xi (z)}\xi '(z)h_0\right)\,d\mu\\
&=B\left[ (\xi , n)|(\xi ', n)\right].
\end{align*}
The compatibility with the inductive limit entails the existence of a sesquilinear extension of $B$ to 
$\mathcal{H}\times\mathcal{H}'$ with
\begin{equation} 
\left| B\left[\xi\, |\,\xi '\right]\right|\leq\sqrt{c}\left\|\xi\right\|_{\mathcal{H}}\left\|\xi '\right\|_{\mathcal{H}'}.
\label{eq2}
\end{equation}
There are some commuting properties between $B$ and $(\tilde{\pi},\tilde{\pi}')$ as follows:
\begin{align*}
B\left[U(\xi , n+1)|U(\xi ',n+1)\right]&=B\left[\left( m_0\left(z^{N^n}\right)\xi (z),n\right)|\left( m_0'\left(z^{N^n}\right)\xi '(z),n\right)\right]\\
&=\int_{\mathbb{T}}R_0^n\left(\overline{m_0\left(z^{N^n}\right)}m_0'\left(z^{N^n}\right)\overline{\xi (z)}\xi '(z)h_0(z)\right)\,d\mu\\
&=\int_{\mathbb{T}}\overline{m_0(z)}m_0'(z)R_0^n\left(\overline{\xi (z)}\xi '(z)h_0(z)\right)\,d\mu\\
&=\int_{\mathbb{T}}R_0\left(R_0^n\left(\overline{\xi (z)}\xi '(z)h_0(z)\right)\right)\,d\mu\\
&=B\left[(\xi , n+1)|(\xi ',n+1)\right].
\end{align*}
So, by density 
\begin{equation}
B\left[U\xi\, |\,U'\xi '\right]=B\left[\xi\, |\, \xi '\right]\, ,\quad (\, \xi\in\mathcal{H}\, , \, \xi '\in\mathcal{H}'\,).
\label{eq3}
\end{equation}
For $f\in\linft$,
\begin{align*}
B\left[\pi (f)(\xi , n)|(\xi ',n)\right]&=B\left[\left(f\left(z^{N^n}\right)\xi(z),n\right)|(\xi ',n)\right]\\
&=\int_{\mathbb{T}}R_0^n\left(\overline{f}\left(z^{N^n}\right)\overline{\xi}(z)\xi '(z)h_0(z)\right)\,d\mu\\
&=B\left[(\xi,n)|\left(\overline{f}\left(z^{N^n}\right)\xi(z),n\right)\right]\\
&=B\left[(\xi,n)|\pi '\left(\overline f\right)(\xi ',n)\right]
\end{align*}
and, again  by density 
\begin{equation}
B\left[\pi (f)\xi |\xi '\right]=B\left[\xi |\pi '\left(\overline f\right)\xi '\right]\, ,\quad (\xi\in\mathcal{H}\, ,\,\xi '\in\mathcal{H}'\,).
\label{eq4}
\end{equation}
As $\varphi=(\openone,0)=\varphi '$ we obtain also
\begin{equation}
B\left[\varphi |\pi '(f)\varphi '\right]=B\left[(\openone ,0)|\pi '(f)(\openone,0)\right]=\int_{\mathbb{T}}f(z)h_0(z)\,d\mu.
\label{eq5}
\end{equation}
\par
Since $B$ is sesquilinear and bounded, there exists $S:\mathcal{H}'\rightarrow\mathcal{H}$ , a bounded linear 
operator with $\left\| S\right\|\leq\sqrt{c}$ such that 
$$B\left[\,\xi\, |\,\xi '\,\right]=\ip{\xi}{S\xi '}\quad ( \, \xi\in\mathcal{H}\, , \, \xi '\in\mathcal{H}'\,).$$
Rewriting (\ref{eq3}) in terms of $S$, one obtains $$SU'=US,$$ 
(\ref{eq4}) gives $$\pi (f)S=S\pi '(f)\quad (f\in\linft),$$ 
 and (\ref{eq5}) yields
$$\ip{\varphi}{\pi (f)S\varphi '}=\int_{\mathbb{T}}fh_0\,d\mu.$$
\par 
For the uniqueness part we will use a lemma which will be useful outside this context too. 
\begin{lemma}
\label{lemma1}
If $n_1, n_2$ are integers and $f_1, f_2\in\linft$ then
\begin{align*}
&\ip{U^{-n_1}\pi\left( f_1\right)\varphi}{S{U'}^{-n_2}\pi '\left(f_2\right)\varphi}=\\
&=\left\{\begin{array}{rcl}
\int_{\mathbb{T}}\overline{f_1}\left(z^{N^{n_2-n_1}}\right)\overline{m}_0^{(n_2-n_1)}(z)f_2(z)h_0(z)\,d\mu &\mbox{ for }&n_2\geq n_1\\
\int_{\mathbb{T}}\overline{f_1}(z){m_0'}^{(n_1-n_2)}(z)f_2\left(z^{N^{n_1-n_2}}\right)h_0(z)\,d\mu &\mbox{ for }&n_1\geq n_2
\end{array}\right.
\end{align*}
\end{lemma}
\begin{proof}
For $n_2\geq n_1$ 
\begin{align*}
\ip{U^{-n_1}\pi\left( f_1\right)\varphi}{S{U'}^{-n_2}\pi '\left(f_2\right)\varphi}&=
\ip{U^{n_2}U^{-n_1}\pi\left( f_1\right)\varphi}{U^{n_2}S{U'}^{-n_2}\pi '\left(f_2\right)\varphi}\\
&=\ip{U^{n_2-n_1}\pi\left( f_1\right)\varphi}{S\pi '\left(f_2\right)\varphi}\\
&=\ip{\pi\left( f_1\left( z^{N^{n_2-n_1}}\right)\right)U^{n_2-n_1}\varphi}{S\pi '\left( f_2\right)\varphi '}\\
\mbox{and using }U^n\varphi=\pi\left(m_0^{(n)}\right)\varphi\\
&=\ip{\pi\left(f_1\left(z^{N^{n_2-n_1}}\right)\right)\pi\left(m_0^{(n_2-n_1)}\right)\varphi}{S\pi '\left(f_2\right)\varphi '}\\
&=\ip{\varphi}{S\pi '\left(\overline{f_1}\left(z^{N^{n_2-n1}}\right)\overline{m}_0^{(n_2-n_1)}(z)f_2(z)\right)\varphi '}\\
&=\int_{\mathbb{T}}\overline{f_1}\left(z^{N^{n_2-n1}}\right)\overline{m}_0^{(n_2-n_1)}(z)f_2(z)h_0(z)\,d\mu.
\end{align*}
For $n_1\geq n_2$ the computation is completely analogous.
\end{proof}
The set 
$$\left\{U^{-n}\pi\left(f\right)\varphi\, |\, n\in\mathbb{N}\, , \, f\in\linft\right\}$$
is dense in $\mathcal{H}$ and similarly for $\mathcal{H}'$, therefore lemma \ref{lemma1} 
implies the uniqueness of $S$. 
\end{proof}
\par
Even more uniqueness holds. In theorem 2.4 in \cite{Jor} it is proved that if $(\pi_1,U_1,\mathcal{H}_1,\varphi_1)$ and 
$(\pi_2,U_2,\mathcal{H}_2,\varphi_2)$ are two cyclic representations of $\mathfrak{A}_N$ corresponding to the same
$h$ then there exists a unique unitary isomorphism $W:\mathcal{H}_1\rightarrow\mathcal{H}_2$ such that
$$W\pi_1(f)=\pi_2(f)W\, ,\quad f\in\linft,\quad WU_1=U_2W, $$
$$W\varphi_1=\varphi_2.$$
\begin{theorem}
\label{th2}
Let $m_0,m_0',h,h',h_0$ be as in theorem \ref{th1}. Suppose $(\pi_i,U_i,\mathcal{H}_i,\varphi_i)$ $i=1,2$ are 
cyclic representations coresponding to $h$, $(\pi_i',U_i',\mathcal{H}_i',\varphi_i')$ $i=1,2$ are cyclic representations 
corresponding to $h'$ and $S_i:\mathcal{H}_i'\rightarrow\mathcal{H}_i$ $i=1,2$ are bounded operators such that
$$S_i\pi_i'(f)=\pi_i(f)S_i,\quad (f\in\linft),\quad S_iU_i'=U_iS_i,\quad (i=1,2)$$
$$\ip{\varphi_i}{S_i\pi_i'(f)\varphi_i'}=\int_{\mathbb{T}}fh_0\,d\mu\quad (i=1,2).$$
Then there exists unique unitary isomorphisms $W:\mathcal{H}_1\rightarrow\mathcal{H}_2$ and $W':\mathcal{H}_1'\rightarrow\mathcal{H}_2'$
such that
$$W\pi_1(f)=\pi_2(f)W,\quad( f\in\linft),\quad WU_1=U_2W,$$
$$W\varphi_1=\varphi_2$$
$$W'\pi_1'(f)=\pi_2'(f)W',\quad (f\in\linft),\quad W'U_1'=U_2'W',$$
$$W\varphi_1'=\varphi_2',$$
$$S_2W'=WS_1.$$
\end{theorem}
\begin{proof}
Let $W,W'$ be given by theorem 2.4 in \cite{Jor} . 
Consider $\tilde S_1=W^{-1}S_2W':\mathcal{H}_1'\rightarrow\mathcal{H}_1$; we prove that $\tilde S_1$ satisfies 
the same conditions as $S_1$ so it must be equal to $S_1$( see theorem \ref{th1}). Let $f\in\linft$. 
\begin{align*}
\tilde S_1\pi_1'(f)&=W^{-1}S_2W'\pi_1'(f)=W^{-1}S_2\pi_2'(f)W'\\
&=W^{-1}\pi_2(f)S_2W'=\pi_1(f)W^{-1}S_2W'\\
&=\pi_1(f)\tilde S_1.
\end{align*}
Similarly $\tilde S_1U_1'=U_1\tilde S_1$.
Also
\begin{align*}
\ip{\varphi_1}{\tilde S_1\pi_1'(f)\varphi_1'}&=
\ip{\varphi_1}{\pi_1(f)W^{-1}S_2W'\varphi_1'}\\
&=\ip{\varphi_1}{\pi_1(f)W^{-1}S_2\varphi_2'}\\
&=\ip{\varphi_1}{W^{-1}\pi_2(f)S_2\varphi_2'}\\
&=\ip{W\varphi_1}{\pi_2(f)S_2\varphi_2'}\\
&=\ip{\varphi_2}{\pi_2(f)S_2\varphi_2'}\\
&=\int_{\mathbb{T}}fh_0\,d\mu.
\end{align*}
\end{proof}
Next we approach the converse problem: that is, given non-sigular $m_0,m_0'\in\linft$ and $h,h'\in\lonet$ 
$h,h'\geq 0$ , $R_{m_0,m_0}h=h$ , $R_{m_0',m_0'}h'=h'$ , consider the cyclic representations 
$(\pi,U,\mathcal{H},\varphi)$ and $(\pi ',U',\mathcal{H}',\varphi ')$. We want to see if an intertwining operator
$S:\mathcal{H}'\rightarrow\mathcal{H}$ induces an eigenvector $h_0\in\lonet$ with 
$$R_{m_0,m_0'}h_0=h_0$$
and
$$\ip{\varphi}{S\pi '(f)\varphi '}=\int_{\mathbb{T}}fh_0\,d\mu,\quad(f\in\linft).$$
The answer is positive and is given in the next theorem.
\begin{theorem}
\label{th3}
Let $m_0,m_0',h,h',(\pi,U,\mathcal{H},\varphi),(\pi ',U',\mathcal{H}',\varphi')$ be as in theorem \ref{th1}.
Suppose $S:\mathcal{H}'\rightarrow\mathcal{H}$ is a bounded operator that satisfies
$$S\pi'(f)=\pi(f)S,\quad (f\in\linft),\quad SU'=US.$$
Then there exists a unique $h_0\in\lonet$ such that 
$$R_{m_0,m_0'}h_0=h_0$$
and
$$\ip{\varphi}{S\pi '(f)\varphi '}=\int_{\mathbb{T}}fh_0\,d\mu,\quad (f\in\linft).$$
Moreover 
$$\left| h_0\right|^2\leq\left\|S\right\|^2hh'\,\mbox{ almost everywhere on }\mathbb{T}.$$ 
\end{theorem}
\begin{proof}
We will need the following result 
\begin{lemma}
If $f_i\in\linft$, $i\in\mathbb{N}$, $f_i$ converges pointwise to $f\in\linft$ $\mu$-a.e. and $\left\|f_i\right\|_{\infty}\leq M<\infty$ for $i\in\mathbb{N}$
then 
$$\pi\left(f_i\right)(\xi)\rightarrow\pi\left(f\right)(\xi)\,\mbox{ in }\mathcal{H}, \quad (\xi\in\mathcal{H}).$$
\end{lemma}
\begin{proof}
Consider first $(\xi,n)\in\mathcal{V}_n$ so $\xi\in\linft$ 
\begin{align*}
\left\|\pi\left(f_i\right)(\xi,n)-\pi\left(f\right)(\xi,n)\right\|_{\mathcal{H}}^2&=
\left\|\left(\left(f_i\left(z^{N^n}\right)-f\left(z^{N^n}\right)\right)\xi(z),n\right)\right\|^2\\
&=\int_{\mathbb{T}}R_{m_0,m_0}^n\left(\left|f_i-f\right|^2\left(z^{N^n}\right)\left|\xi\right|^2(z)h(z)\right)\,d\mu\\
&=\int_{\mathbb{T}}\left|f_i-f\right|^2(z)R_{m_0,m_0}^n\left(\left|\xi\right|^2h\right)\,d\mu\rightarrow 0
\end{align*}
by Lebesque's dominated convergence theorem. 
\par
The set
$$\mathcal{V}=\left\{(\xi,n)\, |\, \xi\in\linft\, ,\, n\in\mathbb{N}\right\}$$
is dense in $\mathcal{H}$.
Fix $\epsilon>0$ and let $a\in\mathcal{H}$. There exists a $b\in\mathcal{V}$ with $\left\|a-b\right\|_{\mathcal{H}}<\epsilon/(3M)$
\begin{align*}
\left\|\pi\left(f_i\right)(a)-\pi\left(f\right)(a)\right\|_{\mathcal{H}}&
\leq\left\|\pi\left(f_i\right)(a)-\pi\left(f_i\right)(b)\right\|_{\mathcal{H}}+
\left\|\pi\left(f_i\right)(b)-\pi\left(f\right)(b)\right\|_{\mathcal{H}}+\\
&+\left\|\pi\left(f\right)(b)-\pi\left(f\right)(a)\right\|_{\mathcal{H}}\\
&\leq\left\|\pi\left(f_i\right)\right\|\left\|a-b\right\|_{\mathcal{H}}+
\left\|\pi\left(f_i\right)(b)-\pi\left(f\right)(b)\right\|_{\mathcal{H}}+\\
&+\left\|\pi\left(f\right)\right\|\left\|a-b\right\|_{\mathcal{H}}
\end{align*}
There is an $n_{\epsilon}$ such that for $i\geq n_{\epsilon}$ one has 
$\left\|\pi\left(f_i\right)(b)-\pi\left(f\right)(b)\right\|_{\mathcal{H}}<\epsilon/3$.
Then for such indices $i$:
$$\left\|\pi\left(f_i\right)(a)-\pi\left(f\right)(a)\right\|_{\mathcal{H}}\leq\left\|f_i\right\|_{\infty}\frac{\epsilon}{3M}
+\frac{\epsilon}{3}+\left\|f\right\|_{\infty}\frac{\epsilon}{3M}\leq\epsilon$$
This concludes the proof of the lemma. 
\end{proof}
\par
Returning to the proof of the theorem, construct the continuous linear functional on $\linft$ by
$$\Lambda:\,f\mapsto\ip{\varphi}{\pi(f)S\varphi'}.$$
Its restriction to the continuous functions on $\mathbb{T}$ is continuous so there is a finite regular Borel measure
 $\nu$ on $\mathbb{T}$ such that 
$$\ip{\varphi}{\pi(f)S\varphi'}=\int_{\mathbb{T}}f\,d\nu\, ,\quad f\in C(\mathbb{T})$$
\par 
Now take $f\in\linft$. Lusin's theorem provides a sequence of continuous functions $f_i$ 
on $\mathbb{T}$ that converges to $f$ pointwise $\mu+|\nu|$-a.e. and with $\left\|f_i\right\|_{\infty}\leq\left\|f\right\|_{\infty}$ for all $i\in\mathbb{N}$.
Using our lemma
$$\ip{\varphi}{\pi(f)S\varphi'}=\lim_{i\to\infty}\ip{\varphi}{\pi(f_i)S\varphi'}=
\lim_{i\to\infty}\int_{\mathbb{T}}f_i\,d\nu=\int_{\mathbb{T}}f\,d\nu,$$
the last equality following from Lebesgue's dominated convergence theorem.
Hence 
$$\ip{\varphi}{\pi(f)S\varphi'}=\int_{\mathbb{T}}f\,d\nu,\quad(f\in\linft).$$
The measure $\nu$ in absolutely continuous, because, if $E$ is a Borel set of $\mu$ measure 
zero, then $\pi\left(\chi_E\right)=0$ so $\nu(E)=\int_{\mathbb{T}}\chi_E\,d\nu=0$. Consequently,
there is an $h_0\in\lonet$ such that $d\nu=h_0\,d\mu$ and we rewrite the previous equation:
$$\ip{\varphi}{\pi(f)S\varphi'}=\int_{\mathbb{T}}fh_0\,d\mu\, , \quad f\in\linft.$$
\par
Next we prove that $R_0h_0=h_0$. Take an arbitrary $f\in\linft$
\begin{align*}
\int_{\mathbb{T}}fh_0\,d\mu &=
\ip{\varphi}{\pi(f)S\varphi'}\\
&=\ip{U\varphi}{U\pi(f)S\varphi'}\\
&=\ip{\pi\left(m_0\right)\varphi}{\pi\left(f\left(z^N\right)\right)SU'\varphi'}\\
&=\ip{\pi\left(m_0\right)\varphi}{S\pi'\left(f\left(z^N\right)\right)\pi'\left(m_0'\right)\varphi'}\\
&=\ip{\pi\left(m_0\right)\varphi}{\pi\left(f\left(z^N\right)m_0'(z)\right)S\varphi'}\\
&=\ip{\varphi}{\pi\left(f\left(z^N\right)m_0'(z)\overline{m_0(z)}\right)S\varphi'}\\
&=\int_{\mathbb{T}}f\left(z^N\right)m_0'(z)\overline{m_0(z)}h_0(z)\,d\mu\\
&=\int_{\mathbb{T}}f(z)R_0h_0\,d\mu.
\end{align*}
As $f$ is arbitrary in $\linft$ this implies that $R_0h_0=h_0$. 
\par
Uniqueness of $h_0$ is clear and we concentrate on last inequality stated in the theorem. 
For all $f,g\in\linft$ we have
\begin{align*}
\left|\int_{\mathbb{T}}\overline fgh_0\,d\mu\right|&=
\left|\ip{\pi(f)\varphi}{\pi(g)S\varphi'}\right|\\
&\leq\left\|\pi\left(f\right)\varphi\right\|_{\mathcal{H}}\left\|S\right\|\left\|\pi'\left(g\right)\varphi'\right\|_{\mathcal{H}}\\
&\leq\left(\int_{\mathbb{T}}\left| f\right|^2h\,d\mu\right)^{\frac 12}\left\|S\right\|
\left(\int_{\mathbb{T}}\left|g\right|^2h'\,d\mu\right)^{\frac 12}
\end{align*}
\par
Since $h_0,h,h'$ are in $\lonet$, almost every point is a Lebesgue density point for all of them. 
Take $z$ a Lebesgue point and $f=g=\chi_I$ where $I$ is a small segment centered at $z$. The inequality above
implies 
$$\left|\int_Ih_0\,d\mu\right|\leq\left\|S\right\|\left(\int_Ih\,d\mu\right)^{\frac 12}
\left(\int_Ih'\,d\mu\right)^{\frac 12},$$
and, dividing by $\mu(I)$, 
$$\left|\frac{1}{\mu(I)}\int_Ih_0\,d\mu\right|\leq\left\|S\right\|\left(\frac{1}{\mu(I)}\int_Ih\,d\mu\right)^{\frac 12}
\left(\frac{1}{\mu(I)}\int_Ih'\,d\mu\right)^{\frac 12}.$$
Then let $I$ shrink to $\left\{z\right\}$ 
$$\left|h_0(z)\right|\leq\left\|S\right\|h^{1/2}(z){h'}^{1/2}(z)$$
and the proof is complete.
\end{proof}
\section{\label{Applications} Applications to wavelets }
\par
We have wavelet representations of $\mathfrak{A}_N$: $ (U,\pi,L^2\left(\mathbb{R}\right),\varphi,m_0)$. Then, by theorem 
\ref{Thmax+b.3}, it is the cyclic representation associated to some positive $h$ with $R_{m_0,m_0}h=h$.
Let's see what the corresponding $h$ is. We must have
$$\ip{\varphi}{\pi(f)\varphi}=\int_{\mathbb{T}}fh\,d\mu,(\quad f\in\linft)$$
and by the unitarity of the Fourier transform 
$$\frac{1}{2\pi}\int_{\mathbb{R}}\overline{\widehat{\varphi}}f\widehat{\varphi}\,dm=\int_{\mathbb{T}}fh\,d\mu$$
and after periodization
$$\int_{\mathbb{T}}f\operatorname*{Per}\left(\left|\widehat{\varphi}\right|^2\right)\,d\mu=\int_{\mathbb{T}}fh\,d\mu$$
which implies $\operatorname*{Per}\left(\left|\widehat{\varphi}\right|^2\right)=h$. Here, for $f\in L^1\left(\mathbb{R}\right)$ 
$$\operatorname*{Per}\left(f\right)(\omega)=\sum_{k\in\mathbb{Z}}f(\omega+2k\pi),\quad \omega\in\left[0,2\pi\right].$$
\par
Hence the wavelet representation is the cyclic representation corresponding to
$\operatorname*{Per}\left(\left|\widehat{\varphi}\right|^2\right)$.
\par
Next, we try to give a neccesary and sufficient condition for the biorthogonality of wavelets.
In the case of biorthogonal wavelets we have the following identity for the filters $m_0,m_0'$ :
$$\frac 1N\sum_{w^N=z}\overline{m_0(w)}m_0'(w)=1,\quad z\in\mathbb{T}$$
which can be rewritten as
$$R_{m_0,m_0'}\openone=\openone.$$
\par
If $\varphi,\varphi'$ are scaling functions corresponding to $m_0,m_0'$ respectively then we see 
that $\operatorname*{Per}\left(\overline{\widehat{\varphi}}\widehat{\varphi '}\right)$ is also an eigenvector 
for $R_{m_0,m_0'}$ .
\par 
Indeed one has the corresponding scaling equations 
$$U\varphi=\pi\left( m_0\right)\varphi\, ,\, U\varphi'=\pi\left(m_0'\right)\varphi'$$
or the Fourier transform versions
$$\widehat\varphi(\omega)=\frac{m_0\left(\frac{\omega}{N}\right)}{\sqrt{N}}\widehat\varphi\left(\frac{\omega}{N}\right)
\, ,\,\widehat\varphi'(\omega)=\frac{m_0'\left(\frac{\omega}{N}\right)}{\sqrt{N}}\widehat\varphi'\left(\frac{\omega}{N}\right).$$
Then
\begin{align*}
\operatorname*{Per}&\left(\overline{\widehat{\varphi}}\widehat{\varphi '}\right)(\omega)
=\sum_{k\in\mathbb{Z}}\overline{\widehat{\varphi}}\widehat{\varphi '}(\omega+2k\pi)\\
&=\sum_{l=0}^{N-1}\sum_{k\in\mathbb{Z}}\overline{\widehat{\varphi}}\widehat{\varphi'}(\omega+2kN\pi+2l\pi)\\
&=\sum_{l=0}^{N-1}\sum_{k\in\mathbb{Z}}\frac{\overline{m_0}\left(\frac{\omega+2l\pi}{N}\right)}{\sqrt{N}}
\overline{\widehat{\varphi}}\left(\frac{\omega+2l\pi}{N}+2k\pi\right)
\frac{m_0'\left(\frac{\omega+2l\pi}{N}\right)}{\sqrt{N}}
\widehat{\varphi'}\left(\frac{\omega+2l\pi}{N}+2k\pi\right)\\
&=\frac 1N\sum_{l=0}^{N-1} \overline{m_0}\left(\frac{\omega+2l\pi}{N}\right)m_0'\left(\frac{\omega+2l\pi}{N}\right)
\sum_{k\in\mathbb{Z}}\overline{\widehat{\varphi}}\widehat{\varphi'}\left(\frac{\omega+2l\pi}{N}+2k\pi\right)\\
&=\frac 1N\sum_{l=0}^{N-1} \overline{m_0}\left(\frac{\omega+2l\pi}{N}\right)m_0'\left(\frac{\omega+2l\pi}{N}\right)
\operatorname*{Per}\left(\overline{\widehat{\varphi}}\widehat{\varphi'}\right)\left(\frac{\omega+2l\pi}{N}\right)\\
&=R_{m_0,m_0'}\left(\operatorname*{Per}\left(\overline{\widehat{\varphi}}\widehat{\varphi'}\right)\right)(\omega).
\end{align*}
\par 
Thus, if we know that $R_{m_0,m_0'}$ has only one eigenvector (up to a multiplicative constant) in some subspace containing $\openone$ and
$\operatorname*{Per}\left(\overline{\widehat{\varphi}}\widehat{\varphi'}\right)$, then we get that 
$$\operatorname*{Per}\left(\overline{\widehat{\varphi}}\widehat{\varphi'}\right)=\openone$$
which is the Fourier transform version of the biorthogonality of $\varphi$ and $\varphi'$. 
\par
We shall see that, under some mild regularity assumptions on $\varphi$ and $\varphi'$, the converse also holds so 
the biorthogonality implies that $R_{m_0,m_0'}$ has a 1-dimensional eigenspace corresponding to the 
eigenvalue $1$. 
\par 
We set up the framework for this converse. Suppose we have cyclic vectors $\varphi$, $\varphi'$ for 
the wavelet representation $\tilde\pi$ of $\mathfrak{A}_N$ on $L^2\left(\mathbb{R}\right)$, satisfying the scaling 
equations
$$U\varphi=\pi\left( m_0\right)\varphi\, ,\, U\varphi'=\pi\left(m_0'\right)\varphi'$$
with $m_0,m_0'$ non-singular in $\linft$ 
\par
The wavelet representation $(U,\pi,L^2\left(\mathbb{R}\right),\varphi,m_0)$ is the cyclic representation 
corresponding to $h=\operatorname*{Per}\left(\left|\widehat{\varphi}\right|^2\right)$. Similarly, for $\varphi'$,
the wavelet representation corresponds to $h'=\operatorname*{Per}\left(\left|\widehat{\varphi}'\right|^2\right)$.

\begin{theorem}
\label{th4}
Let $m_0,m_0',\varphi,\varphi'$ as above. Suppose the following conditions are satisfied:
\begin{enumerate}
\item \label{th4cond1}
$$\widehat{\varphi}(0)\neq 0\neq\widehat{\varphi'}(0)$$
\item \label{th4cond2}
The integer translates of $\varphi$ and $\varphi'$ form Riesz bases for their corresponding linear spans.
\item \label{th4cond3} 
$\widehat{\varphi}$ and $\widehat{\varphi'}$ are continuous at 0 and
$$\sum_{k\neq 0}\left|\widehat{\varphi}\right|^2(x+2k\pi)\rightarrow 0\, \mbox{ as } x\rightarrow 0,$$
$$\sum_{k\neq 0}\left|\widehat{\varphi'}\right|^2(x+2k\pi)\rightarrow 0\, \mbox{ as } x\rightarrow 0.$$
\item \label{th4cond4}
$\varphi$ and $\varphi'$ are biorthogonal or equivalently 
$$\sum_{k\in\mathbb{Z}}\overline{\widehat{\varphi}}(x+2k\pi)\widehat{\varphi'}(x+2k\pi)=1,\quad\mbox{a.e. on }\mathbb{R}.$$
\end{enumerate}
Then there exist exactly one (up to a constant multiple) solution for 
$$R_{m_0,m_0'}h_0=h_0,\quad h_0\in\linft$$
which is continuous at $z=1$ (The solution is $h_0=\openone ).$
\end{theorem}
\begin{proof}
Employing Schwarz' inequality we have 
\begin{align*}
1&=\left|\operatorname*{Per}\left(\overline{\widehat{\varphi}}\widehat{\varphi '}\right)\right|^2(\omega)=
\left|\sum_{k\in\mathbb{Z}}\overline{\widehat{\varphi}}\widehat{\varphi '}(\omega+2k\pi)\right|^2\\
&\leq\left(\sum_{k\in\mathbb{Z}}\left|\widehat{\varphi}\right|^2(\omega+2k\pi)\right)
\left(\sum_{k\in\mathbb{Z}}\left|\widehat{\varphi'}\right|^2(\omega+2k\pi)\right)\\
&=\left(\operatorname*{Per}\left|\widehat{\varphi}\right|^2\right)\left(\operatorname*{Per}\left|\widehat{\varphi'}\right|^2\right)(\omega)
=h(\omega)h'(\omega).
\end{align*}
Suppose $h_0$ is a solution in $\linft$ of $R_{m_0,m_0'}h_0=h_0$ which is continuous at $z=1$. Since 
$h_0\in\linft$, there is a $c<\infty$ such that $\left|h_0\right|^2\leq c$ almost everywhere on $\mathbb{T}$ 
and the previous inequality implies 
$$\left|h_0\right|^2\leq chh'.$$
Therefore, by theorem \ref{th1}, $h_0$ induces an operator $S:L^2\left(\mathbb{R}\right)\rightarrow
L^2\left(\mathbb{R}\right)$ such that 
$$US=SU\, , \, \pi(f)S=S\pi(f),\quad(f\in\linft),$$
$$\ip{\varphi}{S\pi(f)\varphi'}=\int_{\mathbb{T}}fh_0\,d\mu.$$
Fourier transform the last equation and then periodize to obtain for $f\in\linft$

$$\int_{\mathbb{T}}fh_0\,d\mu=\frac{1}{2\pi}\ip{\widehat{\varphi}}{f\widehat{\left(S\varphi'\right)}}
=\int_{\mathbb{T}}\operatorname*{Per}\left(\overline{\widehat{\varphi}}\widehat{\left(S\varphi'\right)}\right)f\,d\mu.$$
So
\begin{equation} \label{eq1th4}
h_0=\operatorname*{Per}\left(\overline{\widehat{\varphi}}\widehat{\left(S\varphi'\right)}\right)
\end{equation}
Also the commuting properties of $S$ imply that $S\varphi'$ is a solution for the scaling equation
$$US\varphi'=\pi\left(m_0'\right)S\varphi'$$
\par
Assume $S\varphi'\neq c\varphi'$ where $c$ is some constant. 
We want to prove that $\widehat{S\varphi'}$ can't be continuous at $0$. Otherwise, consider the Fourier version of
the scaling equation 
$$\widehat{S\varphi'}(\omega)=\frac{m_0'\left(\frac{\omega}{N}\right)}{\sqrt{N}}\widehat{S\varphi'}\left(\frac{\omega}{N}\right)$$
and by induction
$$\widehat{S\varphi'}(\omega)=\left[\prod_{i=1}^n\frac{m_0'\left(\frac{\omega}{N^i}\right)}{\sqrt{N}}\right]\widehat{S\varphi'}\left(\frac{\omega}{N^n}\right).$$
 A similar equation can be constructed for $\widehat{\varphi'}$. If $\widehat{S\varphi'}$ is continuous at 0 then
$$\frac{\widehat{S\varphi'}(x)}{\widehat{\varphi'}(x)}=\lim_{n\to\infty}\frac{\widehat{S\varphi'}\left(\frac{x}{N^n}\right)}{\widehat{\varphi'}\left(\frac{x}{N^n}\right)}
=\frac{\widehat{S\varphi'}(0)}{\widehat{\varphi'}(0)}$$
So $S\varphi'=c\varphi$ with $c=\frac{\widehat{S\varphi'}(0)}{\widehat{\varphi'}(0)}$, a contradiction.
\par 
On the other hand, from (\ref{eq1th4}) we get
\begin{equation}\label{eq2th4}
\overline{\widehat{\varphi}}(x)\widehat{S\varphi'}(x)=h_0(x)-\sum_{k\neq 0}\overline{\widehat{\varphi}}(x+2k\pi)
\widehat{S\varphi'}(x+2k\pi).
\end{equation}
$\overline{\widehat{\varphi}}$ and $h_0$ are continuous at $x=0$. 
We prove that the sum in (\ref{eq2th4}) converges to $0$ as $x\to 0$. By the Schwarz inequality
\begin{equation}\label{eq3th4}
\left|\sum_{k\neq 0}\overline{\widehat{\varphi}}\widehat{S\varphi '}(x+2k\pi)\right|\leq
\left(\sum_{k\neq 0}\left|\widehat{\varphi}\right|^2(x+2k\pi)\right)
\left(\sum_{k\neq 0}\left|\widehat{S\varphi'}\right|^2(x+2k\pi)\right).
\end{equation}
We try to bound the second factor. Since $S$ commutes with $\pi$ and $U$, the same is true for 
$S^*$ and $S^*S$. By theorem \ref{th3}, $S^*S$ induces some $h_0'$ with $R_{m_0',m_0'}h_0'=h_0'$,
$\left|h_0'\right|^2\leq c{h'}^2$ for some $c>0$ and
$$\ip{\varphi'}{\pi(f)S^*S\varphi'}=\int_{\mathbb{T}}fh_0'\,d\mu,(\quad f\in\linft).$$
Then 
$$\ip{S\varphi'}{\pi(f)S\varphi'}=\int_{\mathbb{T}}fh_0'\,d\mu.$$

Using again the Fourier transform and periodization we obtain
$$\operatorname*{Per}\left|\widehat{S\varphi'}\right|^2=h_0'\leq\sqrt{c}h'.$$
Since $\varphi'$ generates a Riesz basis, (see \cite{Dau92}) there is a $B<\infty$ such that 
$$h'=\operatorname*{Per}\left|\widehat{\varphi'}\right|^2\leq B$$
Thus $\operatorname*{Per}\left|\widehat{S\varphi'}\right|^2$ is bounded and, using the hypothesis (\ref{th4cond3}),
 in (\ref{eq3th4}) we obtain that 
$$\lim_{x\to 0}\left|\sum_{k\neq 0}\overline{\widehat{\varphi}}\widehat{S\varphi '}(x+2k\pi)\right|=0.$$
Now apply this to (\ref{eq2th4}) and use the fact that $\widehat{\varphi}$ is continuous at 0 with $\widehat{\varphi}(0)\neq 0$
to conclude that $\widehat{S\varphi'}$ is continuous 0, again a contradiction which leads to 
$S\varphi'=c\varphi'$. Without loss of generality, take the constant $c$ to be 1. 
$$\int_{\mathbb{T}}fh_0\,d\mu=\ip{\varphi}{\pi(f)S\varphi'}=\ip{\varphi}{\pi(f)\varphi'}=\int_{\mathbb{T}}f\,d\mu$$
the last equality follows from (\ref{th4cond4}) using the usual Fourier-periodization technique.
The equality holds for all $f\in\linft$ so $h_0=\openone$ .
\end{proof}
\begin{corollary}\label{remark1}
If $\varphi$ and $\varphi'$ are compactly supported, biorthogonal, 
$$m_0(0)=\sqrt{N}=m_0'(0)\, , m_0\left(\frac{2k\pi}{N}\right)=0=m_0'\left(\frac{2k\pi}{N}\right)\, ,\quad k\in\left\{1,...,N-1\right\}$$
then $\openone$ is the only solution of $R_{m_0,m_0'}h_0=h_0$ which is continuous at $z=1$ (up to a multiplicative constant).
\end{corollary}
\begin{proof}
We have to check the conditions of theorem \ref{th4}.  
The Fourier coefficients of $\operatorname*{Per}\left|\hat\varphi\right|^2$ are given by
$$\int_{\mathbb{T}}e_k\operatorname*{Per}\left|\hat\varphi\right|^2\,d\mu=\ip{\varphi}{\pi\left(e_k\right)\varphi}=
\ip{\varphi}{\varphi(\omega -k)},$$
where $e_k=e^{-ik\theta}$. Therefore the coefficients are zero except for a finite number of them so
$\operatorname*{Per}\left|\hat\varphi\right|^2$ is a trigonometric polynomial. The same is true for 
$\operatorname*{Per}\left|\hat{\varphi'}\right|^2$. Also the fact that $\varphi$ and $\varphi'$ are compactly 
supported implies that $\hat{\varphi},\hat{\varphi'}$ are continuous. From the scaling equation we obtain that 
$\widehat{\varphi(0)}=1=\widehat{\varphi'}(0)$. 
\par
For $k\in\mathbb{Z}$, $k\neq0$, we can write $k=N^pl$ with $l=Nq+r$, $r\in\left\{1,...,N-1\right\}$. Then
\begin{align*}
\hat{\varphi}(2k\pi)&=\hat{\varphi}\left(2N^pl\pi\right)\\
&=\frac{m_0\left(2N^{p-1}l\pi\right)}{\sqrt{N}}\ldots\frac{m_0\left(2Nl\pi\right)}{\sqrt{N}}
\frac{m_0(2l\pi)}{\sqrt{N}}\frac{m_0\left(\frac{2l\pi}{N}\right)}{\sqrt{N}}\hat{\varphi}\left(\frac{2l\pi}{N}\right)=0
\end{align*}
because 
$$m_0\left(\frac{2l\pi}{N}\right)=m_0\left(\frac{2(Nq+r)\pi}{N}\right)=m_0\left(\frac{2r\pi}{N}\right)=0.$$
\par 
Thus $\hat{\varphi}(2k\pi)=0$ for $k\neq0$. This shows that $\operatorname*{Per}\left|\hat\varphi\right|^2(0)=1$. Then
$$\sum_{k\neq0}\left|\hat{\varphi}\right|^2(x+2k\pi)=\operatorname*{Per}\left|\hat\varphi\right|^2(x)-\hat{\varphi}(x)\rightarrow
\operatorname*{Per}\left|\hat\varphi\right|^2(0)-\hat{\varphi}(0)=0\,\mbox{ as } x\rightarrow 0.$$
The same argument applies to $\hat{\varphi'}.$ 
\par
Condition (\ref{th4cond2}) is obtain as follows: look at the first inequality in the proof of theorem \ref{th4}. 
We have 
$$1\leq\left(\operatorname*{Per}\left|\widehat{\varphi}\right|^2\right)\left(\operatorname*{Per}\left|\widehat{\varphi'}\right|^2\right)(\omega)$$
As both factors are trigonometric polynomials they are bounded by a common constant $0<A<\infty$ . Hence 
$$1/A\leq\left(\operatorname*{Per}\left|\widehat{\varphi}\right|^2\right)\leq A$$ 
which implies that the translates of $\varphi$ form a Riesz basis for their linear span (see \cite{Dau92}). 
Similarly for $\varphi'$.
\end{proof}

\section{\label{Examples} Some examples }
 We know that in the case of a quadrature mirror filter $m_0$ for which the transfer operator 
$R_{m_0,m_0}$ has 1 as a simple eigenvalue in $C(\mathbb{T})$, the cyclic representation that 
corresponds to the constant function $\openone$ is in fact the wavelet representation on $L^2(\mathbb{R})$.
Then the commutant of this representation is in one-to-one correspondence with $\linft$-solutions 
for $R_{m_0,m_0}h=h$. We will describe this commutant and give the form of all corresponding 
$\linft$-solutions. 
\par
 We recall that the wavelet representation of $\mathfrak{A}_N$ is generated by
$$U :\psi\mapsto\frac{1}{\sqrt{N}}\psi\left(\frac{x}{N}\right)\quad ,\quad 
V:\psi\mapsto\psi(x-1)$$
$V=\pi(z)$ where $\pi$ is the representation of $\linft$ given by
$$\left(\pi(f)\psi\right)^{\widehat{ }}=f\widehat\psi ,\quad(f\in\linft , \psi\in\ltwor).$$ 
 It will be useful to consider the representation in Fourier space and in this case $U$ has 
the form
\begin{equation}\label{eqexu}
\widehat U\widehat{\psi}(x)=\sqrt{N}\widehat{\psi}(Nx).
\end{equation}
 Thus we have another representation of $\mathfrak{A}_N$ on $L^2(\mathbb{R})$ , which is 
generated by $\widehat U$ given in (\ref{eqexu}) and the representation $\widehat{\pi}$ of $\linft$ given by
$$\widehat{\pi}(f)\widehat{\psi}=f\widehat{\psi},\quad (f\in\linft\, ,\, \widehat{\psi}\in L^2(\mathbb{R})).$$
This representation is equivalent to the wavelet representation via the Fourier transform. 
\begin{theorem}\label{th41}
The commutant of $(\widehat{\pi},\widehat U,L^2(\mathbb{T}))$ is 
$$\left\{ M_f\,|\, f\in L^{\infty}(\mathbb{R})\, ,\, f(Nx)=f(x) \mbox{ a.e. on } \mathbb{R}\right\}$$
where $M_f(\psi)=f\psi$ for all $\psi\in L^{2}(\mathbb{R}), f\in\linfr.$
\end{theorem}
\begin{proof} 
(The theorem is valid even in a more general case , see \cite{Li00}) 
 Let $A$ be an operator that commutes with $\widehat U$ and $\widehat{\pi}(f)$ for all $f\in\linft$. We prove first 
that $A$ commutes with $M_g$, where $g\in L^{\infty}(\mathbb{R})$ is periodic of period 
$2N\pi$. 
\par
Indeed, let $f(x)=g(Nx)$ for $x\in\mathbb{R}$. Then $f$ is $2\pi$-periodic and bounded so 
$\widehat{\pi}(f)$ commutes 
with $A$. Then $A$ commutes also with $\widehat U^{-1}\widehat{\pi}(f)\widehat U$. But
$$(\widehat U^{-1}\widehat{\pi}(f)\widehat U\psi)(x)=\widehat U^{-1}\left(f(x)\sqrt{N}\psi(Nx)\right)=f\left(\frac{x}{N}\right)\psi(x)
=g(x)\psi(x).$$
 So $\widehat U^{-1}\widehat{\pi}(f)\widehat U=M_g$. It follows by induction that $A$ commutes with any multiplication by a $2\pi N^k$-periodic, 
function in $L^{\infty}(\mathbb{R})$.
 Now, consider $f\in L^{\infty}(\mathbb{R})$. We claim that $A$ commutes with $M_f$. Define 
$f_n(x)=f(x)$ on $[-\pi N^n,\pi N^n]$ and extended it on $\mathbb{R}$ by $2\pi N^n$-periodicity. 
First we prove that $M_{f_n}$ converges to $M_f$ in the strong operator topology. For this, 
take $\psi\in L^{2}(\mathbb{R})$. 
Then
\begin{align*}
\left\| M_{f_n}\psi -M_f\psi\right\|^2&=\int_{\mathbb{R}}\left|f_n-f\right|^2\left|\psi\right|^2\,dx\\
&=\int_{\left| x\right|\geq\pi N^n}\left|f_n-f\right|^2\left|\psi\right|^2\,dx\\
&\leq\left( 2\left\|f\right\|_{\infty}\right)^2\int_{\mathbb{R}}\chi_{\left\{\left| x\right|\geq\pi N^n\right\}}\left|\psi\right|^2\,dx\\
&\rightarrow 0\,\mbox{ as } n\rightarrow\infty.
\end{align*}
\par 
Thus $M_f$ is the limit of $M_{f_n}$ in the strong operator topology and consequently $A$ will commute 
with $M_f$ also. Using then theorem IX.6.6 in \cite{Con90} we obtain that $A=M_f$ for some $f\in 
L^{\infty}\left(\mathbb{R}\right)$. Then, the fact that $A$ and $\widehat U$ commute implies:
$$f(Nx)=f(x)\,\mbox{ a.e. on }\mathbb{R}.$$ 
\par 
This proves one inclusion, the other one is a straightforward verification. 
\end{proof}
\par
Using this theorem we can find all solutions to $R_{m_0,m_0'}h=h$ as follows:
\begin{corollary}\label{cor42}
Suppose we have the wavelet representations $(U,\pi,L^{2}\left(\mathbb{R}\right),\varphi,m_0)$ and
$(U,\pi,L^{2}\left(\mathbb{R}\right),\varphi',m_0')$. Let $h=\operatorname*{Per}\left(\left|\widehat{\varphi}\right|^2\right)$ 
and $h'=\operatorname*{Per}\left(\left|\widehat{\varphi}'\right|^2\right)$. Then each solution $h_0\in
L^{1}\left(\mathbb{R}\right)$ for $R_{m_0,m_0'}h_0=h_0$ with $\left|h_0\right|^2\leq chh'$ for some positive constant $c$, 
has the form 
$$h_0=\operatorname*{Per}\left(f\overline{\widehat{\varphi}}\widehat{\varphi'}\right)$$
for some $f\in L^{\infty}\left(\mathbb{R}\right)$ with $f(Nx)=f(x)$ a.e. 
\par
Conversely, any such $h_0$ is an $L^{1}\left(\mathbb{R}\right)$-solution for $R_{m_0,m_0'}h_0=h_0$ 
and $\left|h_0\right|^2\leq chh'$ for some $c>0$. 
\end{corollary}
\begin{proof}
\par
The cyclic representations corresponding to $h$ and $h'$ are the wavelet representations given 
in the hypothesis. Hence the intertwining operators are in fact the ones in the commutant of 
the wavelet representation. 
We will transfer everything into the Fourier space by applying the Fourier transform and then use 
theorem \ref{th41} and the results in section \ref{Main}. 
\par
If $h_0$ is as given then there is an operator $A$ in the commutant of the representation 
such that 
$$\ip{\varphi}{\pi(g)A\varphi'}=\int_{\mathbb{T}}gh_0\,d\mu,\quad (g\in\linft).$$
But after the Fourier transform, we saw that $\widehat A$ has the form $\widehat{A}=M_f$ where 
$f\in L^{\infty}\left(\mathbb{R}\right)$ with $f(Nx)=f(x)$ a.e. . Therefore 
$$\frac{1}{2\pi}\ip{\widehat{\varphi}}{gf\widehat{\varphi'}}=\int_{\mathbb{T}}gh_0\,d\mu,\quad (g\in\linft)$$
and after periodization we obtain that
$$h_0=\operatorname*{Per}\left(f\overline{\widehat{\varphi}}\widehat{\varphi'}\right).$$
\par
Conversely, it is easy to see that, when $f$ is given, 
$$\frac{1}{2\pi}\ip{\widehat{\varphi}}{gM_f\widehat{\varphi'}}=\int_{\mathbb{T}}gh_0\,d\mu,\quad (g\in\linft)$$
and as $M_f$ is in the commutant, the rest follows.
\end{proof}
\begin{example}
 In the sequel we consider $N=2$, $m_0(z)=\frac{1}{\sqrt{2}}\left(1+z^p\right)$, $p$ being an odd integer. 
This example appears also in \cite{BraJo}. 
We try to find out the solutions for $R_{m_0,m_0}h=h$. 
\par
 It is easy to see that 
$$R_{m_0,m_0}\openone=\openone$$
Also if $\varphi=\frac{1}{p}\chi_{(0,p)}$ then  
$$U\varphi=\pi\left(m_0\right)\varphi.$$
Then the wavelet representation $(U,\pi,L^2\left(\mathbb{T}\right),\varphi,m_0)$ is the cyclic 
representation corresponding to $h_{\varphi}=\operatorname*{Per}\left(\left|\widehat{\varphi}\right|^2\right)$,
$$\widehat{\varphi}(x)=e^{-ip\frac{x}{2}}\frac{\sin\frac{px}{2}}{\frac{px}{2}}.$$
\par
Using the identity
\begin{equation}\label{eq42}
\sum_{n\in\mathbb{Z}}\frac{1}{(t+2\pi n)^2}=\frac{1}{4\sin^2\left(\frac{t}{2}\right)},\quad(t\in\mathbb{R}).
\end{equation}
we obtain
$$h_{\varphi}(t)=\frac{1}{p^2}\frac{\sin^2\left(\frac{pt}{2}\right)}{\sin^2\left(\frac{t}{2}\right)},\quad(t\in\mathbb{R}).$$
\par
We try to construct the cyclic representation corresponding to $\openone$. Let $\rho=e^{i\frac{2\pi}{p}}$ and 
for $\eta\in\mathbb{T}$ define $\alpha_{\eta}(f)(z)=f(\eta z)$ for all $f\in\linft$ and $z\in\mathbb{T}$. The essential 
observation is that we have the following identity:
\begin{equation}\label{eq43}
\sum_{k=0}^{p-1}\alpha_{\rho^k}\left(h_{\varphi}\right)=1\, ,\mbox{ on }\mathbb{T}
\end{equation}
This identity follows from the following computation:
\begin{align*}
\sum_{k=0}^{p-1}\alpha_{\rho^k}\left(h_{\varphi}\right)(t)&=
\sum_{k=0}^{p-1}h_{\varphi}\left(t+\frac{2k\pi}{p}\right)\\
&=\frac{4}{p^2}\sin^2\left(\frac{pt}{2}\right)\sum_{k=0}^{p-1}\sum_{n\in\mathbb{Z}}\frac{1}{\left(t+\frac{2k\pi}{p}+2n\pi\right)^2}\\
&=4\sin^2\left(\frac{pt}{2}\right)\sum_{l\in\mathbb{Z}}\frac{1}{\left(pt+2\pi l\right)^2}\,\quad(l=k+pn)\\
&=1\, , \mbox{(using (\ref{eq42})).}
\end{align*}
\par
We construct now the cyclic representation that corresponds to $\openone$. 
Let $$\mathcal{H}_1=\underbrace{\ltwor\oplus\ltwor\oplus\ldots\oplus\ltwor}_{p \mbox{ times}}$$
For $f\in\linft$ define 
$$\pi_1(f)\left(\xi_0,\ldots,\xi_{p-1}\right)=\left(\pi(\alpha_{\rho^0}(f))(\xi_0),\pi(\alpha_{\rho^1}(f))(\xi_1),\ldots,
\pi(\alpha_{\rho^{p-1}}(f))(\xi_{p-1})\right),$$
$$U_1\left(\xi_0,\ldots,\xi_{p-1}\right)=\left(U\xi_{\sigma(0)},U\xi_{\sigma(1)},\ldots,U\xi_{\sigma(p-1)}\right),$$
where $\pi$ and $U$ come from the wavelet representation on $\ltwor$ and $\sigma$ is the permutation of the set 
$\left\{0,\ldots,p-1\right\}$ given by $\rho^{2k}=\rho^{\sigma(k)}$ so $\sigma(k)=2k\mbox{ mod }p$. 
\par
 $\pi_1$ is a representation of $\linft$, $U_1$ is unitary, and a short computation, based on the fact that 
$\alpha_{\rho^k}\left(f\left(z^2\right)\right)=\left(\alpha_{\rho^{\sigma(k)}}(f)\right)\left(z^2\right)$, shows that   
$$U_1\pi_1(f)=\pi_1\left(f\left(z^2\right)\right)U_1,\quad (f\in\linft)$$
\par
Define $\varphi_1=(\varphi,\varphi,\ldots,\varphi)$. Since $\alpha_{\rho^k}(m_0)=m_0$ for all $k\in\left\{0,\ldots,p-1\right\}$
and $U\varphi=\pi\left(m_0\right)\varphi$,
$$U_1\varphi_1=\pi_1\left(m_0\right)\varphi_1.$$ And, finally, we check that $\openone$ is the eigenfunction that induces this representation.
As $\ip{\varphi}{\pi(f)\varphi}=\int_{\mathbb{T}}fh_{\varphi}\,d\mu$, we have 
\begin{align*}
\ip{\varphi_1}{\pi_1(f)\varphi_1}&=\sum_{k=0}^{p-1}\int_{\mathbb{T}}\alpha_{\rho^k}(f)h_{\varphi}\,d\mu\\
&=\sum_{k=0}^{p-1}\int_{\mathbb{T}}f\alpha_{\rho^{-k}}(h_{\varphi})\,d\mu\\
&=\int_{\mathbb{T}}f\sum_{l=0}^{p-1}\alpha_{\rho^l}(h_{\varphi})\,d\mu\\
&=\int_{\mathbb{T}}f\openone\,d\mu\, ,\quad ( \mbox{ by (\ref{eq43}) } )
\end{align*}
\par
Now we Fourier transform everything and try to find the commutant of the representation. After 
the Fourier transform, $\widehat{\pi}_1$ and $\widehat{U}_1$ have the same form as before, only now, $\widehat{\pi}(f)$ is the multiplication by 
$f$, ($f\in\linft$), and 
$$\widehat{U}\psi(x)=\sqrt{2}\psi(2x),\quad (\psi\in\ltwor)$$
\par
So consider $A:\mathcal{H}_1\rightarrow\mathcal{H}_1$ that commutes with the representation,
$A=\left(A_{ij}\right)_{i,j=0}^{p-1}$. 
\par
Consider $g\in\linft$ of period $\frac{2\pi}{p}$. Then $\alpha_{\rho^k}(g)=g$ for all $k$. Hence
$$\widehat{\pi}_1(g)(\xi_0,\ldots,\xi_{p-1})=(\widehat{\pi}(g)(\xi_0),\ldots,\widehat{\pi}(g)(\xi_{p-1})).$$
Also, since $\sigma$ is permutation, there is an $M$ such that $\sigma^M(k)=k$ for all $k$ so
$$\widehat{U}_1^M(\xi_0,\ldots,\xi_{p-1})=(\widehat{U}^M\xi_0,\ldots,\widehat{U}^M\xi_{p-1}).$$
\par
Note that $P_i$ , the projection on the $i$-th component, commutes with $\widehat{\pi}_1(g)$ and $\widehat{U}_1^M$. 
Then $P_iAP_j$ commutes with $\widehat{\pi}_1(g)$ and $\widehat{U}_1^M$, and this implies that $A_{ij}$ commutes with 
$\widehat{\pi}(g)$ and $\widehat{U}^M$. Repeating the argument in theorem \ref{th41}, we obtain that $A_{ij}=M_{f_{ij}}$ 
for some $f_{ij}\in\linfr$.
\par
Now take $f\in\linft$ arbitrary. The fact that $\widehat{\pi}_1(f)$ and $A$ commute can be rewritten as
$$\sum_{j=0}^{p-1}f_{ij}\alpha_{\rho^j}(f)\xi_j=\alpha_{\rho^i}(f)\sum_{j=0}^{p-1}f_{ij}\xi_j,\quad (i\in\left\{0,\ldots,p-1\right\},(\xi_j)\in\mathcal{H}_1)$$ 
Fix $k$ and let $\xi_j=0$ for $j\neq k$. Then
$$f_{ik}\alpha_{\rho^k}(f)\xi_k=\alpha_{\rho^i}(f)f_{ik}\xi_k$$
and this implies that $f_{ik}=0$ for $i\neq k$.
\par
Now the fact that $A$ commutes with $\widehat{U}_1$ implies
$$f_{ii}(x)=f_{\sigma(i)\sigma(i)}(2x),\mbox{ a.e. on }\mathbb{R}\, ,\,(i\in\left\{0,\ldots,p-1\right\})$$
\par
A simple check shows that any $A$ of this form commutes with the representation, therefore the commutant is given by :
\begin{align*}
\widehat{\pi}_1'&=\left\{A:\mathcal{H}_1\rightarrow\mathcal{H}_1\,|\,A(\xi_0,\ldots,\xi_{p-1})=(f_0\xi_0,\ldots,f_{p-1}\xi_{p-1}),\right.\\
&\left. f_i\in\linfr , f_i(x)=f_{\sigma(i)}(2x) \mbox{ a.e .}\right\}.
\end{align*}
\par
Using this and section \ref{Main} we obtain that the $\linft$-solutions for $R_{m_0,m_0}h=h$ must satisfy the 
identity
$$\frac{1}{2\pi}\sum_{k=0}^{p-1}\ip{\widehat{\varphi}}{\alpha_{\rho^k}(f)f_k\widehat{\varphi}}=\int_{\mathbb{T}}fh,\quad (f\in\linft).$$
for some $f_k\in\linfr$ with $f_k(x)=f_{\sigma(k)}(2x)$ a.e. on $\mathbb{R}$. So 
$$h=\sum_{k=0}^{p-1}\operatorname*{Per}\left(f_k\left(x+\frac{2k\pi}{p}\right)\left|\widehat{\varphi}\right|^2\left(x+\frac{2k\pi}{p}\right)\right).$$
\par
To conclude our analysis we try to find the continuous eigenfunctions. So let's take $h$ continuous, having the 
form mentioned above. We want to see what conclusions can be drawn on $f_i$. We prove that they are 
constants. 
\par
Fix $k\in\left\{0,\ldots,p-1\right\}$.
\begin{align*}
h(x)&=f_k\left(x+\frac{2k\pi}{p}\right)\left|\widehat{\varphi}\right|^2\left(x+\frac{2k\pi}{p}\right)+\\
&+\sum_{l\in\mathbb{Z}\setminus\left\{0\right\}}f_k\left(x+\frac{2k\pi}{p}+2l\pi\right)\left|\widehat{\varphi}\right|^2\left(x+\frac{2k\pi}{p}+2l\pi\right)+\\
&+\sum_{j\in\left\{0,\ldots,p-1\right\}\setminus\left\{k\right\}}\sum_{l\in\mathbb{Z}}f_j\left(x+\frac{2j\pi}{p}+2l\pi\right)\left|\widehat{\varphi}\right|^2\left(x+\frac{2j\pi}{p}+2l\pi\right).
\end{align*}
Denote the first sum by $A(x)$ and the second one by $B(x)$. 
$$\left|A(x)\right|\leq M\sum_{l\in\mathbb{Z}\setminus\left\{0\right\}}\left|\widehat{\varphi}\right|^2\left(x+\frac{2k\pi}{p}+2l\pi\right),$$
where $M=\sup_i\sup_x\left|f_i(x)\right|$.
This last sum is a continuous function, as it is the difference between $h_{\varphi}\left(x+\frac{2k\pi}{p}\right)$ 
and $\left|\widehat{\varphi}\right|^2\left(x+\frac{2k\pi}{p}\right)$, and its value at $-\frac{2k\pi}{p}$ 
is 0, since
$$\left|\widehat{\varphi}\right|^2\left(-\frac{2k\pi}{p}+\frac{2k\pi}{p}+2l\pi\right)=0,\quad l\neq0.$$
Similarly,
$$\left|B(x)\right|\leq M\sum_{j\in\left\{0,\ldots,p-1\right\}\setminus\left\{k\right\}}\operatorname*{Per}\left|\widehat{\varphi}\right|^2\left(x+\frac{2j\pi}{p}\right)$$
and also this function is continuous and $0$ at $x=-\frac{2k\pi}{p}$ ( look at $h_{\varphi}$ ).
Thus 
$$A(x)+B(x)\rightarrow 0,\,\mbox{ as } x\rightarrow -\frac{2k\pi}{p}.$$
Since $h$ is continuous, the following identity holds 
$$\lim_{x\rightarrow -\frac{2k\pi}{p}}f_k\left(x+\frac{2k\pi}{p}\right)\left|\widehat{\varphi}\right|^2\left(x+\frac{2k\pi}{p}\right)=h\left(-\frac{2k\pi}{p}\right).$$
But $\left|\widehat{\varphi}\right|^2$ is also continuous and $\left|\widehat{\varphi}\right|^2(0)=1$ so 
$$\lim_{x\rightarrow 0}f_k(x)=h\left(-\frac{2k\pi}{p}\right).$$
\par
On the other hand since $\sigma$ is a permutation, there is a $K$ such that $\sigma^K(i)=i$ for all $i\in\left\{0,\ldots,p-1\right\}$. 
Then by induction $f_k(x)=f_{\sigma^K(k)}(2^Kx)=f_k(2^Kx)$ a.e. on $\mathbb{R}$. This, coupled with the 
limit at 0, makes $f_k$ constant $a_k$. Then
$$h=\sum_{k=0}^{p-1}a_k\operatorname*{Per}\left|\widehat{\varphi}\right|^2\left(x+\frac{2k\pi}{p}\right)$$
so 
$$h=\sum_{k=0}^{p-1}a_k\alpha_{\rho^k}\left(h_{\varphi}\right)$$
with $a_k=a_{\sigma(k)}$ for all $k\in\left\{0,\ldots,p-1\right\}$.
\par
We want to give a basis for the space of the continuous eigenfunctions $h$. 
For this, note that $a_k$ is constant for $k$ in a cycle of $\sigma$. So let 
$O_1,\ldots,O_{c(p)}$ be the cycles of $\sigma$. Then each continuous $h$ will have the form 
$$h=\sum_{k=1}^{c(p)}b_k\sum_{l\in O_k}\alpha_{\rho^l}\left(h_{\varphi}\right).$$ 
The functions $\sum_{l\in O_k}\alpha_{\rho^l}\left(h_{\varphi}\right)$ can be seen to be linearly independent 
if we observe that the set of zeroes are $\left\{\rho^{-l}\,|\,l\in \left\{0,\ldots,p-1\right\}\setminus O_k\right\}$. 
\par
Also, it is easy to see that the cyclic representation associated to $$h_{O_k}=\sum_{l\in O_k}\alpha_{\rho^l}\left(h_{\varphi}\right)$$ 
is given by $P_{O_k}\pi_1$, $P_{O_k}U_1$, where $P_{O_k}$ is the projection on the components in $O_k$. 
\par
Take now , for example $p=9$. So $m_0(z)=\frac{1}{\sqrt{2}}\left(1+z^9\right)$, $\varphi=\frac{1}{9}\chi_{(0,9)}$. 
$$h_{\varphi}(t)=\frac{1}{9^2}\frac{\sin^2\left(\frac{9t}{2}\right)}{\sin^2\left(\frac{t}{2}\right)}.$$
it induces the wavelet representation on $\ltwor$. 
\par
$$\sigma=\left(\begin{array}{ccccccccc}
0&1&2&3&4&5&6&7&8\\
0&2&4&6&8&1&3&5&7
\end{array}\right)
$$
The cycles are $O_1=\{0\}$, $O_2=\{1,2,4,5,7,8\}$, $O_3=\{3,6\}$. 
\par
$h_{O_1}=h_{\varphi}$ of course.\\
Observe that 
\begin{align*}
h_{O_1}(x)+h_{O_3}(x)&=\frac{4}{9^2}\sum_{k=0}^{2}\sum_{n\in\mathbb{Z}}\frac{\sin^2\left(\frac{9\left(x+\frac{2k\pi}{3}+2n\pi\right)}{2}\right)}{\left(x+\frac{2k\pi}{3}+2n\pi\right)^2}\\
&=\frac{4}{3^2}\sin^2\left(\frac{9x}{2}\right)\sum_{k=0}^{2}\sum_{n\in\mathbb{Z}}\frac{1}{\left(3x+2\pi(k+3n)\right)^2}\\
&=\frac{4}{3^2}\sin^2\left(\frac{9x}{2}\right)\sum_{l\in\mathbb{Z}}\frac{1}{\left(3x+2\pi l\right)^2}\\
&=\frac{1}{3^2}\frac{\sin^2\left(\frac{9x}{2}\right)}{\sin^2\left(\frac{3x}{2}\right)}.
\end{align*}
Also 
$$h_{O_1}+h_{O_2}+h_{O_3}=\openone.$$
 Therefore a basis for the continuous eigenfunctions is 
$$\left\{\openone , \frac{1}{3^2}\frac{\sin^2\left(\frac{9x}{2}\right)}{\sin^2\left(\frac{3x}{2}\right)} , \frac{1}{9^2}\frac{\sin^2\left(\frac{9x}{2}\right)}{\sin^2\left(\frac{x}{2}\right)}\right\}.$$

\end{example}

\begin{acknowledgements}
The author wants to thank professor Palle E.T. Jorgensen for his suggestions and his 
always kind and prompt help.
\end{acknowledgements}

\end{document}